\documentclass[12pt, reqno]{amsart}
\usepackage{amsmath, amsthm, amscd, amsfonts, amssymb, graphicx, color}
\usepackage[bookmarksnumbered, colorlinks, plainpages]{hyperref}

\newtheorem{theorem}{Theorem}[section]
\newtheorem{lemma}[theorem]{Lemma}
\newtheorem{proposition}[theorem]{Proposition}
\newtheorem{corollary}[theorem]{Corollary}
\theoremstyle{definition}
\newtheorem{definition}[theorem]{Definition}
\newtheorem{example}[theorem]{Example}

\theoremstyle{remark}
\newtheorem{remark}[theorem]{Remark}
\numberwithin{equation}{section}

\allowdisplaybreaks
\begin{document}

\title[Continuous Riesz bases in Hilbert $C^*$-modules]
{Continuous Riesz bases in Hilbert $C^*$-modules}

\author[H.Ghasemi]{Hadi Ghasemi }
\address{Hadi Ghasemi \\ Department of Mathematics and Computer
Sciences, Hakim Sabzevari University, Sabzevar, P.O. Box 397, IRAN}
\email{ \rm h.ghasemi@hsu.ac.ir}
\author[T.L. Shateri]{Tayebe Lal Shateri }
\address{Tayebe Lal Shateri \\ Department of Mathematics and Computer
Sciences, Hakim Sabzevari University, Sabzevar, P.O. Box 397, IRAN}
\email{ \rm  t.shateri@hsu.ac.ir; shateri@ualberta.ca}
\thanks{*The corresponding author:
t.shateri@hsu.ac.ir; shateri@ualberta.ca (Tayebe Lal Shateri)}
 \maketitle

\begin{abstract}
The paper is devoted to continuous frames and Riesz bases in Hilbert $C^*$-modules. we define a continuous Riesz basis for Hilbert $C^*$-modules and give some results about them.
\vskip 3mm
\end{abstract}

\section{Introduction}
Frame theory is nowadays a fundamental research area in mathematics, computer science and engineering with many interesting applications in a variety of different fields. Frames were first introduced by Duffin and Schaeffer \cite{DS} in the context of nonharmonic fourier series. Then Daubecheies, Grassman and Mayer \cite{DG} reintroduced and developed them. The concept of a generalization of frames to a family indexed by some locally compact space endowed with a Radon measure was proposed by G. Kaiser \cite{KAI} and independently by Ali, Antoine and Gazeau \cite{AAG}. These frames are known as continuous frames. For a discussion of continuous frames, we refer to Refs.\cite{RND, RDD}. Arefijamaal and et al. \cite{AKR} introduced continuous Riesz bases and give some equivalent conditions for a continuous frame to be a continuous Riesz basis.

One reason to study frames in Hilbert $C^*$-modules is that there are some differences between
Hilbert spaces and Hilbert $C^*$-modules. For example, in general, every bounded operator on a Hilbert space has an unique adjoint, while this fact not hold for bounded operators on a Hilbert $C^*$-module. 
Thus it is more difficult to make a discussion of the theory of Hilbert $C^*$-modules than that of Hilbert spaces in general. We refer the readers to \cite{LAN}, for more details on Hilbert $C^*$-modules. Frank and Larson \cite{FL} presented a general approach to the frame theory in Hilbert $C^*$-modules. The theory of frames has been extended from Hilbert spaces to Hilbert $C^*$-modules, see \cite{AR,JI,FL,RR,SH}.
 
The paper is organized as follows. First, we recall the basic definitions and some notations about Hilbert $C^*$-modules, and we also give some properties of them. Also, we recall the notion of continuous frames in Hilbert $C^*$-modules and their operators. In section 3, we define continuous Riesz bases in Hilbert $C^*$-modules and we give some results aboute them.
\section{Preliminaries}
First, we recall some definitions and basic properties of Hilbert $C^*$-modules. We give only a brief introduction to the theory of Hilbert $C^*$-modules to make our explanations self-contained. For comprehensive accounts, we refer to  \cite{LAN,OLS}. Throughout this paper, $\mathcal A$ shows a unital $C^*$-algebra.
\begin{definition}
A \textit{pre-Hilbert module} over unital $C^*$-algebra $\mathcal A$ is a complex vector space $U$ which is also a left $\mathcal A$-module equipped with an $\mathcal A$-valued inner product $\langle .,.\rangle :U\times U\to \mathcal A$ which is $\mathbb C$-linear and $\mathcal A$-linear in its first variable and satisfies the following conditions:\\
$(i)\; \langle f,f\rangle \geq 0$,\\
$(ii)\; \langle f,f\rangle =0$  iff $f=0$,\\
$(iii)\; \langle f,g\rangle ^*=\langle g,f\rangle ,$\\
$(iv)\; \langle af,g\rangle=a\langle f,g\rangle ,$\\
for all $f,g\in U$ and $a\in\mathcal A$.
\end{definition}
A pre-Hilbert $\mathcal A$-module $U$ is called \textit{Hilbert $\mathcal A$-module} if $U$ is complete with respect to the topology determined by the norm $\|f\|=\|\langle f,f\rangle \|^{\frac{1}{2}}$.

By \cite[Example 2.46]{JI}, if $\mathcal A$ is a $C^*$-algebra, then it is a Hilbert $\mathcal A$-module with respect to the inner product
$$\langle a,b\rangle =ab^*,\quad (a,b\in \mathcal A).$$
\begin{example}
\cite[Page 237]{OLS} Let $l^2(\mathcal A)$ be the set of all sequences $\{a_n\}_{n\in \mathbb N}$ of elements of a $C^*$-algebra $\mathcal A$ such that the series $\sum_{n=1}^{\infty}a_na_n^*$ is convergent in $\mathcal A$. Then $l^2(\mathcal A)$ is a Hilbert $\mathcal A$-module with respect to the pointwise operations and inner product defined by
\begin{equation*}
\langle \{a_n\}_{n\in \mathbb N},\{b_n\}_{n\in \mathbb N}\rangle =\sum_{n=1}^{\infty}a_nb_n^*.
\end{equation*}
\end{example}
In the following lemma the \textit{Cauchy-Schwartz inequality} reconstructed in Hilbert $C^*$-modules.
\begin{lemma}
\cite[Lemma 15.1.3]{OLS} (\textbf{Cauchy-Schwartz inequality}) Let $U$ be a Hilbert $C^*$-modules over a unital $C^*$-algebra $\mathcal A$. Then
\begin{equation*}
\|\langle f,g\rangle \|^{2}\leq\|\langle f,f\rangle \|\;\|\langle g,g\rangle \|,
\end{equation*}
for all $f,g\in U$.
\end{lemma}
\begin{definition}
\cite[Page 8]{LAN} Let $U$ and $V$ be two Hilbert $C^*$-modules over a unital $C^*$-algebra $\mathcal A$. A map $T:U\to V$ is said to be \textit{adjointable} if there exists a map $T^{*}:V\to U$ satisfying
$$\langle Tf,g\rangle =\langle f,T^*g\rangle $$,
for all $f\in U, g\in V$. Such a map $T^*$ is called the \textit{adjoint} of $T$. By $End_{\mathcal A}^*(U)$ we denote the set of all adjointable maps on $U$.
\end{definition}
It is surprising that an adjointable operator is automatically linear and bounded.
\begin{lemma} \label{UKR}
\cite[Lemma 1.1]{XS} Let $U$ and $V$ be two Hilbert $C^*$-modules over a unital $C^*$-algebra $\mathcal A$ and $T\in End_{\mathcal A}^*(U,V)$ has closed range. Then $T^*$ has closed range and 
$$U=Ker(T)\oplus R(T^{*})\;\;\;\;\; ,\;\;\;\;\; V=Ker(T^{*})\oplus R(T)$$.
\end{lemma}
\begin{lemma}\label{ISIN}
\cite[Lemma 0.1]{AD} Let $U$ and $V$ be two Hilbert $C^*$-modules over a unital $C^*$-algebra $\mathcal A$ and $T\in End_{\mathcal A}^*(U,V)$. Then\\
(i)\;If $T$ is injective and $T$ has closed range, then the adjointable map $T^{*}T$ is invertible and 
\begin{equation*}
\|(T^{*}T)^{-1}\|^{-1}\leq T^{*}T\leq\| T\|^{2}.
\end{equation*}
(i)\;If $T$ is surjective, then the adjointable map $TT^{*}$ is invertible and 
\begin{equation*}
\|(TT^{*})^{-1}\|^{-1}\leq TT^{*}\leq\| T\|^{2}.
\end{equation*}
\end{lemma}
Now, we introduce continuous frames in Hilbert $C^*$-modules over a unital $C^*$-algebra $\mathcal A$, and then we give some results for these frames.  

Let $\mathcal Y$ be a Banach space, $(\mathcal X,\mu)$ a measure space, and $f:\mathcal X\to \mathcal Y$ a measurable function. The integral of the Banach-valued function $f$ has been defined by Bochner and others. Most properties of this integral are similar to those of the integral of real-valued functions (see \cite{DAN,YOS}). Since every $C^*$-algebra and Hilbert $C^*$-module is a Banach space, hence we can use this integral in these spaces. In the following, we assume that $\mathcal A$ is a unital $C^*$-algebra and $U$ is a Hilbert $C^*$-module over $\mathcal A$ and $(\Omega ,\mu)$ is a measure space. 
\begin{definition}
\cite{GHSH1,ROS} Let $(\Omega ,\mu)$ be a measure space and $\mathcal A$ is a unital $C^*$-algebra. We define,
\begin{equation*}
L^{2}(\Omega ,A)=\lbrace\varphi :\Omega \to A\quad ;\quad \Vert\int_{\Omega}\vert(\varphi(\omega))^{*}\vert^{2} d\mu(\omega)\Vert<\infty\rbrace .
\end{equation*}
For any $\varphi ,\psi \in L^{2}(\Omega , A)$, the inner product is defined by $\langle \varphi ,\psi\rangle = \int_{\Omega}\langle\varphi(\omega),\psi(\omega)\rangle d\mu(\omega)$ and the norm is defined by $\|\varphi\|=\|\langle \varphi,\varphi\rangle \|^{\frac{1}{2}}$. It was shown in \cite{LAN} $L^{2}(\Omega , A)$ is a Hilbert $\mathcal A$-module.
\end{definition}
Continuous frames for Hilbert $\mathcal A$-modules are defined as follows.
\begin{definition}
\cite{GHSH1,ROS} A mapping $F:\Omega \to U$ is called a continuous frame for $U$ if\\
$(i)\; F$ is weakly-measurable, i.e, for any $f\in U$, the mapping $\omega\longmapsto\langle f,F(\omega)\rangle $ is measurable on $\Omega$.\\
$(ii)$ There exist constants $A,B>0$ such that
\begin{equation}\label{eq1}
A\langle f,f\rangle \leq \int_{\Omega}\langle f,F(\omega)\rangle \langle F(\omega),f\rangle d\mu(\omega)\leq B\langle f,f\rangle  ,\quad (f\in U). 
\end{equation}
\end{definition}

The constants $A,B$ are called \textit{lower} and \textit{upper} frame bounds, respectively. The mapping $F$ is called \textit{Bessel} if the right inequality in \eqref{eq1} holds and is called \textit{tight} if  $A=B$.
\begin{definition}
\cite{GHSH} A continuous frame $F:\Omega \to U$ is called \textit{exact} if for every measurable subset $\Omega_{1}\subseteq\Omega$ with
$0<\mu(\Omega_{1})<\infty$, the mapping $F:\Omega\backslash\Omega_{1} \to U$ is not a continuous frame for $U$.
\end{definition}
\begin{example}
Let
$U=\Big\{
\begin{pmatrix}
a&0&0\\
0&0&b
\end{pmatrix}
: a,b\in \mathbb C\Big\}$, and $\mathcal A=\Big\{
\begin{pmatrix}
x&0\\
0&y
\end{pmatrix}
: x,y\in \mathbb C\Big\}$ which is a $C^*$-algebra. We define the inner product
\begin{equation*}
 \begin{array}{ll}
\langle .,.\rangle:U\times U\;\to \quad \mathcal A \\
\qquad\; (M,N)\longmapsto M(\overline{N})^t.
 \end{array}
 \end{equation*}
This inner product makes $U$ a $C^*$-module on $\mathcal A$. We consider a measure space $(\Omega=[0,1],\mu)$ whose $\mu$ is the Lebesgue measure. Also $F:\Omega\to U$ defined by 
$F(\omega)=\begin{pmatrix}
\sqrt{3}\omega&0&0\\
0&0&\sqrt{3}\omega
\end{pmatrix}$, for any $\omega\in \Omega$.\\
For each $f=
\begin{pmatrix}
a&0&0\\
0&0&b
\end{pmatrix}
\in U$, we have
\begin{align*}
\int_{[0,1]}\langle f,F(\omega)\rangle\langle F(\omega),f\rangle d\mu(\omega)&=\int_{[0,1]}
3\omega ^2\begin{pmatrix}
|a|^2&0\\
0&|b|^2
\end{pmatrix}d\mu(\omega)\\
&=\begin{pmatrix}
|a|^2&0\\
0&|b|^2
\end{pmatrix}
=\langle f,f\rangle .
\end{align*}
Therefore $F$ is a continuous tight frame with bounds $A=B=1$.
\end{example}
The following operators for Bessel mappings and continuous frames in Hilbert $C^{\ast}$-modules are defined in \cite{GHSH}.
\begin{definition}
Let $F:\Omega \to U$ be a Bessel mapping . Then \\
(i)\;The \textit{synthesis operator} or \textit{pre-frame operator} $T_{F}:L^{2}(\Omega , A)\;\to U$ weakly defined by
\begin{equation}
\langle T_{F}\varphi ,f\rangle =\int_{\Omega}\varphi(\omega)\langle F(\omega),f\rangle d\mu(\omega),\quad (f\in U).
\end{equation}
(ii)\; The adjoit of $T$, called The \textit{analysis operator} $T^{\ast}_{F}:U\;\to L^{2}(\Omega , A)$ is defined by
\begin{equation}
(T^{\ast}_{F}f)(\omega)=\langle f ,F(\omega)\rangle ,\quad (\omega\in \Omega).
\end{equation}
\end{definition}
The pre-frame operator is a well defined, surjective, adjointable $\mathcal A$-linear map and is bounded with $\| T\|\leq\sqrt{B}$ and the analysis operator $T^{\ast}_{F}:U\;\to L^{2}(\Omega , A)$ is injective and has closed range \cite{GHSH}.
\begin{definition}
Let $F:\Omega \to U$ be a continuous frame for Hilbert $C^{\ast}$-module $U$. Then the frame operator $S_{F}:U\;\to U$ is weakly defined by
\begin{equation}
\langle S_{F}f ,f\rangle =\int_{\Omega}\langle f,F(\omega)\rangle\langle F(\omega),f\rangle d\mu(\omega),\quad (f\in U).
\end{equation}
\end{definition}
In \cite{GHSH} prove that the frame operator $S_F$ is positive, adjointable, selfadjoit and invertible and the lower and  the upper bounds of $F$ are respectively $\Vert S^{-1}\Vert^{-1}$ and $\Vert T\Vert^{2}$.
Now we introduce the concept of the duals of continuous frames in Hilbert $C^{\ast}$-modules and give some important properties of continuous frames and their duals.
\begin{definition}
\cite{GHSH} Let $F:\Omega \to U$ be a continuous Bessel mapping. A continuous Bessel mapping $G:\Omega \to U$ is called a \textit{dual} for $F$ if
\begin{equation*}
f= \int_{\Omega}\langle f,G(\omega)\rangle F(\omega)d\mu(\omega) ,\;\;\;\;\;(f\in U),
\end{equation*}
or
\begin{equation}\label{eq2}
\langle f,g\rangle= \int_{\Omega}\langle f,G(\omega)\rangle\langle F(\omega),g\rangle d\mu(\omega) ,\;\;\;\;\;(f,g\in U).
\end{equation}
In this case $(F,G)$ is called a \textit{dual pair}. If $T_{F}$ and $T_{G}$ denote the synthesis operators of $F$ and $G$, respectively, then \eqref{eq2} is equivalent to $T_{F}T^{*}_{G}=I_{U}$. The condition 
\begin{equation*}
\langle f,g\rangle= \int_{\Omega}\langle f,G(\omega)\rangle\langle F(\omega),g\rangle d\mu(\omega),\;\;\;\;\;(f,g\in U),
\end{equation*}
is equivalent
\begin{equation*}
\langle f,g\rangle= \int_{\Omega}\langle f,F(\omega)\rangle\langle G(\omega),g\rangle d\mu(\omega),\;\;\;\;\;(f,g\in U),
\end{equation*}
because $T_{F}T^{*}_{G}=I_{U}$ if and only if $T_{G}T^{*}_{F}=I_{U}$.\\
Also, by reconstructin formula we have
\begin{equation*}
f=S^{-1}Sf=S^{-1}\int_{\Omega}\langle f,F(\omega)\rangle F(\omega) d\mu(\omega)=\int_{\Omega}\langle f,F(\omega)\rangle S^{-1}F(\omega) d\mu(\omega),
\end{equation*}
and
\begin{equation*}
f=SS^{-1}f=\int_{\Omega}\langle S^{-1}f,F(\omega)\rangle F(\omega) d\mu(\omega)=\int_{\Omega}\langle f,S^{-1}F(\omega)\rangle F(\omega) d\mu(\omega).
\end{equation*}
Then $S^{-1}F$ is a dual for $F$, which is called \textit{canonical dual}.
\end{definition}
 \section{Continuous Riesz bases in Hilbert $C^*$-modules}
In this section, we introduce the concept of continuous Riesz bases in  Hilbert $C^{\ast}$-modules and give some important properties of them. First, we give the notion of a Riesz-type frame that is introduced in \cite{GHSH}.
\begin{definition}
Let $F:\Omega \to U$ be a continuous frame for Hilbert $C^{\ast}$-module $U$. If $F$ has only one dual, we call $F$ a \textit{Riesz-type frame}.
\end{definition}
\begin{theorem}\cite[Theorem 3.4]{GHSH}
Let $F:\Omega \to U$ be a continuous frame for Hilbert $C^{\ast}$-module $U$ over a unital $C^*$-algebra $\mathcal A$. Then $F$ is a Riesz-type frame if and only if the analysis operator $T^{*}_{F}:U\to L^{2}(\Omega ,A)$ is onto.
\end{theorem}
\begin{definition}
Let $U$ be a Hilbert $C^{\ast}$-module $U$ over a unital $C^*$-algebra $\mathcal A$. A Bessel mapping $F:\Omega \to U$ is called a \textit{$\mu$-complete} if
$$\lbrace f\in U\;\; ;\;\;\langle f,F(\omega)\rangle =0\;\; a.e.\;[\mu]\rbrace =\lbrace 0\rbrace$$.
\end{definition}
Now, we define a continuous Riesz basis for Hilbert $C^{\ast}$-modules.
\begin{definition}
Let $U$ be a Hilbert $C^{\ast}$-module $U$ over a unital $C^*$-algebra $\mathcal A$. A mapping $F:\Omega \to U$ is called a \textit{continuous Riesz basis} for Hilbert $C^{\ast}$-module $U$, if the following conditions are satisfied:\\
$(i)\; F$ is weakly-measurable, i.e, for any $f\in U$, the mapping $\omega\longmapsto\langle f,F(\omega)\rangle $ is measurable on $\Omega$.\\
$(ii)\;$ $F$ is $\mu$-complete.\\
$(iii)\;$ There are two constants $A,B>0$ such that
\begin{equation}
A\|\int_{\Omega_{1}}\vert \varphi(\omega)^{*}\vert^{2} d\mu(\omega)\|^{\dfrac{1}{2}} \leq \|\int_{\Omega_{1}}\varphi(\omega)F(\omega) d\mu(\omega)\|\leq B\|\int_{\Omega_{1}}\vert \varphi(\omega)^{*}\vert^{2}d\mu(\omega)\|^{\dfrac{1}{2}},
\end{equation}
for every $\varphi\in L^{2}(\Omega , A)$ and measurable subset $\Omega_{1}\subseteq\Omega$ with  $\mu(\Omega_{1})<+\infty$.
\end{definition}
\begin{remark}
Let $F:\Omega \to U$ be a continuous Riesz basis for Hilbert $C^{\ast}$-module $U$. Define 
\begin{align*}
T:L^{2}(\Omega & ,A)\longrightarrow U\\
& \varphi\longmapsto \int_{\Omega}\varphi(\omega)F(\omega) d\mu(\omega)
\end{align*}
Then $T$ is well-defined, adjointable map with $T^{*}f=\lbrace\langle f,F(\omega)\rangle\rbrace_{\omega\in\Omega}$ and bounded such that
\begin{equation*}
 A\Vert\varphi\Vert^{2}\leq\Vert T\varphi\Vert^{2}\leq B\Vert\varphi\Vert^{2}.
\end{equation*}
Hence $F$ is a continuous Bessel mapping. Also by $\mu$-completeness of $F$ we have
\begin{equation*}
 Ker(T^{*})=\lbrace f\in U \;\; ;\;\;\langle f,F(\omega)\rangle =0 \;\;\; \forall\omega\in\Omega\rbrace =\lbrace 0\rbrace ,
\end{equation*}
so by lemma \ref{UKR}, $ R(T)=Ker(T^{*})^{\perp}=U$. Then $T$ is onto and by \cite[theorem 2.15]{GHSH}, $F$ is a continuous frame for $U$.
\end{remark}
\begin{definition}
Let $U$ be a Hilbert $C^{\ast}$-module $U$ over a unital $C^*$-algebra $\mathcal A$. A Bessel mapping $F:\Omega \to U$ is said to be \textit{$L^{2}$-independent} if for $\varphi\in L^{2}(\Omega ,A)$,\\
$\int_{\Omega}\varphi(\omega)F(\omega) d\mu(\omega)=0$ implies that $\varphi (\omega) =0$, for each $\omega\in\Omega$.
\end{definition}
We give the following result.
\begin{theorem}
Let $F:\Omega \to U$ be a continuous frame for Hilbert $C^{\ast}$-module $U$ over a unital $C^*$-algebra $\mathcal A$ with bounds $A,B>0$. Then the following are equivalent:\\
$(i)\; F$ is a continuous Riesz basis.\\
$(ii)\; F$ is $\mu$-complete and $L^{2}$-independent.
\end{theorem}
\begin{proof}
$(i)\;\Longrightarrow\;(ii)$\; Let $F$ be a continuous Riesz basis and $\int_{\Omega}\varphi(\omega)F(\omega) d\mu(\omega)=0$ for some $\varphi\in L^{2}(\Omega ,A)$. Since 
\begin{equation*}
 A\Vert\int_{\Omega}\vert \varphi(\omega)^{*}\vert^{2} d\mu(\omega)\Vert\leq\Vert\int_{\Omega}\varphi(\omega)F(\omega) d\mu(\omega)\Vert^{2}=0,
\end{equation*}
so 
\begin{equation*}
\langle\lbrace\varphi(\omega)\rbrace_{\omega\in\Omega} ,\lbrace\varphi(\omega)\rbrace_{\omega\in\Omega}\rangle =\int_{\Omega}\vert \varphi(\omega)^{*}\vert^{2} d\mu(\omega)=0.
\end{equation*}
Hence $\lbrace\varphi(\omega)\rbrace_{\omega\in\Omega}=0$ and $\varphi =0$ i.e. $F$ is $L^{2}$-independent.

$(ii)\;\Longrightarrow\;(i)$ Let $F$ be a $L^{2}$-independent continuous frame for Hilbert $C^{\ast}$-module $U$ with bounds $A,B>0$. For $\varphi\in L^{2}(\Omega ,A)$ and measurable subset $\Omega_{1}\subseteq\Omega$ with $\mu(\Omega_{1})<+\infty$, put $f=\int_{\Omega_{1}}\varphi(\omega)F(\omega) d\mu(\omega)$. Then we have,
\begin{equation*}
f=\int_{\Omega_{1}}\varphi(\omega)F(\omega) d\mu(\omega)=\int_{\Omega}\varphi(\omega)\chi_{\Omega_{1}}(\omega) F(\omega) d\mu(\omega).
\end{equation*}
Also $f=\int_{\Omega}\langle f,S^{-1}F(\omega)\rangle F(\omega) d\mu(\omega)$ where $S$ is the continuous frame operator of $F$.\\
Since $F$ is $L^{2}$-independent, so
\begin{equation*}
\varphi(\omega)\chi_{\Omega_{1}}(\omega) = \langle f,S^{-1}F(\omega)\rangle ,\;\;\;\;\;(\omega\in\Omega).
\end{equation*}
and by \cite[corollary 2.11]{GHSH},
\begin{equation*}
B^{-1}\langle f,f\rangle\leq\langle S^{-1}f,f\rangle\leq A^{-1}\langle f,f\rangle ,
\end{equation*}
and so
\begin{equation*}
A\Vert\langle S^{-1}f,f\rangle\Vert\leq\Vert\langle f,f\rangle\Vert\leq B\Vert\langle S^{-1}f,f\rangle\Vert .
\end{equation*}
Now we show that $\langle S^{-1}f,f\rangle =\int_{\Omega_{1}}\vert\varphi(\omega)^{*}\vert^{2} d\mu(\omega)$.
\begin{align*}
\langle S^{-1}f,f\rangle & =\langle f,S^{-1}f\rangle\\
& =\int_{\Omega_{1}}\langle f,\varphi(\omega)S^{-1}F(\omega)\rangle d\mu(\omega)\\
& =\int_{\Omega_{1}}\langle f,S^{-1}F(\omega)\rangle\varphi(\omega)^{*} d\mu(\omega)\\
& =\int_{\Omega_{1}}\varphi(\omega)\varphi(\omega)^{*} d\mu(\omega) =\int_{\Omega_{1}}\vert\varphi(\omega)^{*}\vert^{2} d\mu(\omega).
\end{align*}
Therefore,
\begin{equation*}
A\|\int_{\Omega_{1}}\vert \varphi(\omega)^{*}\vert^{2} d\mu(\omega)\| \leq \|\int_{\Omega_{1}}\varphi(\omega)F(\omega) d\mu(\omega)\|^{2}\leq B\|\int_{\Omega_{1}}\vert \varphi(\omega)^{*}\vert^{2}d\mu(\omega)\|,
\end{equation*}
i.e. $F$ is a continuous Riesz basis for $U$ with bounds $A,B$.
\end{proof}  
\begin{theorem}\label{exact}
Let $F:\Omega \to U$ be a continuous frame for Hilbert $C^{\ast}$-module $U$ over a unital $C^*$-algebra $\mathcal A$. If $F$ is a continuous Riesz basis for $U$, then it is a continuous exact frame.
\end{theorem}
\begin{proof}
Let $\Omega_{1}\subseteq\Omega$ be a measurable subset of $\Omega$ with $0<\mu(\Omega_{1})<\infty$. For $\varphi =\chi_{\Omega_{1}}\in L^{2}(\Omega ,A)$ we have,
\begin{align*}
\Vert\int_{\Omega_{1}}F(\omega) d\mu(\omega)\Vert^{2} & =\Vert\int_{\Omega_{1}}\chi_{\Omega_{1}}(\omega)F(\omega) d\mu(\omega)\Vert^{2} \\
& \geq A\Vert\int_{\Omega_{1}}\vert\chi_{\Omega_{1}}(\omega)\vert^{2} d\mu(\omega)\Vert\\
& =A\Vert\mu(\Omega_{1})\Vert >0.
\end{align*}
Hence $\int_{\Omega_{1}}F(\omega) d\mu(\omega)\neq 0$.

Now suppose that $F:\Omega\setminus\Omega_{1} \to U$ is a continuous frame for $U$. Then by completeness of $F\mid_{\Omega\setminus\Omega_{1}}$ there exists $\varphi_{0}\in L^{2}(\Omega\setminus\Omega_{1} ,A)$ such that
$$\int_{\Omega_{1}}F(\omega) d\mu(\omega)=\int_{\Omega\setminus\Omega_{1}}\varphi_{0}(\omega)F(\omega) d\mu(\omega).$$
Define $\varphi:\Omega \to A$ where
\begin{equation*}
\varphi(\omega)=
\begin{cases}
\varphi_{0}(\omega) & \omega\in\Omega\setminus\Omega_{1}\\
1 & \omega\in\Omega_{1}.
\end{cases}
\end{equation*}
Then $\varphi\in L^{2}(\Omega ,A)$ and
\begin{equation*}
\int_{\Omega}\chi_{\Omega_{1}}(\omega)F(\omega) d\mu(\omega)=\int_{\Omega}\varphi(\omega)F(\omega) d\mu(\omega),
\end{equation*}
so $\int_{\Omega}(\varphi -\chi_{\Omega_{1}})(\omega)F(\omega) d\mu(\omega) =0$. Hence $L^{2}$-independent shows that $\varphi =\chi_{\Omega_{1}}$ and so $\varphi_{0} =0$.

Therefore
\begin{equation*}
\int_{\Omega_{1}}F(\omega) d\mu(\omega)=\int_{\Omega\setminus\Omega_{1}}\varphi_{0}(\omega)F(\omega) d\mu(\omega)=0,
\end{equation*}
which is a contradiction.
\end{proof}
\begin{proposition}\label{Vmap}
Let $F:\Omega \to U$ be a continuous Bessel mapping for Hilbert $C^{\ast}$-module $U$ over a unital $C^*$-algebra $\mathcal A$ with pre-frame operator $T$. Suppose that $F$ is $\mu$-complete and the mapping
\begin{align*}
V:\; L^{2} & (\Omega ,A)\longrightarrow L^{2}(\Omega ,A)\\
& \varphi\longmapsto\int_{\Omega}\varphi(\omega)\langle F(\omega),F(.)\rangle d\mu(\omega)
\end{align*}
defines a bounded, adjointable and invertible operator. Then $F$ is a continuous Riesz basis for $U$.
\end{proposition}
\begin{proof}
Since $F$ is Bessel, so the synthesis operator $T$ is well-defined and bounded and adjointable with $T^{*}f=\lbrace\langle f,F(\omega)\rangle\rbrace_{\omega\in\Omega}$ for $f\in U$.\\
Also $T^{*}T=V$, because
\begin{align*}
(T^{*}T)(\varphi) & =T^{*}(\int_{\Omega}\varphi(\omega) F(\omega) d\mu(\omega))\\
& = \lbrace\langle\int_{\Omega}\varphi(\omega) F(\omega) d\mu(\omega),F(\gamma)\rangle\rbrace_{\gamma\in\Omega}\\
& = \lbrace\int_{\Omega}\varphi(\omega)\langle F(\omega) ,F(\gamma)\rangle d\mu(\omega)\rbrace_{\gamma\in\Omega}.
\end{align*}
Since $T$ is bounded, so there exist $B>0$ such that $\Vert T\varphi\Vert^{2}\leq B\Vert\varphi\Vert^{2}$ i.e.
\begin{equation*}
\Vert\int_{\Omega}\varphi(\omega) F(\omega) d\mu(\omega)\Vert^{2}\leq B\Vert\int_{\Omega}\vert\varphi(\omega)^{*}\vert^{2} d\mu(\omega)\Vert.
\end{equation*}
Since $T^{*}T$ is positive, so
\begin{align*}
\Vert\int_{\Omega}\varphi(\omega) F(\omega) d\mu(\omega)\Vert^{2} & =\Vert T\varphi\Vert^{2}\\
& = \Vert\langle T^{*}T\varphi ,\varphi\rangle\Vert\\
& = \Vert\langle (T^{*}T)^{\dfrac{1}{2}}\varphi ,(T^{*}T)^{\dfrac{1}{2}}\varphi\rangle\Vert\\
& = \Vert (T^{*}T)^{\dfrac{1}{2}}\varphi\Vert^{2}\geq \Vert (T^{*}T)^{\dfrac{-1}{2}}\Vert^{-2}\Vert\varphi\Vert^{2}.
\end{align*}
Therefore $F$ is continuous Riesz basis with lower and upper bounds $\Vert (T^{*}T)^{\dfrac{-1}{2}}\Vert^{-2}$,$B$ respectively.
\end{proof}
 \begin{theorem}\label{R-type}
Let $F:\Omega \to U$ be a continuous frame for Hilbert $C^{\ast}$-module $U$ over a unital $C^*$-algebra $\mathcal A$ with pre-frame operator $T$. Then $F$ is a continuous Riesz basis for $U$ if and only if $F$ is a Riesz-type frame.
\end{theorem}
\begin{proof}
$(\Longrightarrow)$\;Let $F_{1}\neq F_{2}$ are two duals of $F$. Then for each $f\in U$,
\begin{align*}
\int_{\Omega}\langle f,F_{1}(\omega)-F_{2}(\omega)\rangle F(\omega) d\mu(\omega) & =\int_{\Omega}\langle f,F_{1}(\omega)\rangle F(\omega) d\mu(\omega) -\int_{\Omega}\langle f,F_{2}(\omega)\rangle F(\omega) d\mu(\omega) \\
& =f-f=0.
\end{align*}
Since $F$ is continuous Riesz basis, so is $L^{2}$-independent and
\begin{equation*}
\langle f,F_{1}(\omega)-F_{2}(\omega)\rangle =0 \;\;\;\Longrightarrow\;\;\;\langle f,F_{1}(\omega)\rangle =\langle f,F_{2}(\omega)\rangle\;\;\;(\omega\in\Omega).
\end{equation*}
Therefore $F_{1}= F_{2}$.

$(\Longleftarrow)$\;Let $F$ be Riesz-type frame and $\varphi\in L^{2}(\Omega ,A)$ such that $\int_{\Omega}\varphi(\omega)F(\omega) d\mu(\omega)=0$.\\
Since $F$ is Riesz-type, so $R(T^{*})=L^{2}(\Omega ,A)$. Also $L^{2}(\Omega ,A)= Ker(T)\oplus R(T^{*})$.
Then
$$\varphi\in Ker(T) = R(T^{*})^{\perp} =\lbrace 0\rbrace ,$$
then $\varphi =0$ and so $F$ is $L^{2}$-independent. Therefore $F$ is a continuous Riesz basis.
\end{proof}
\begin{corollary}
Let $F:\Omega \to U$ be a continuous frame for Hilbert $C^{\ast}$-module $U$ over a unital $C^*$-algebra $\mathcal A$. If $F$ is a Riesz-type frame, then it is a continuons exact frame.
\end{corollary}
Duo to the Theorem \ref{R-type}, the converse of the Proposition \ref{Vmap} holds as follows.
\begin{corollary}
Let $F:\Omega \to U$ be a continuous Riesz basis for Hilbert $C^{\ast}$-module $U$ over a unital $C^*$-algebra $\mathcal A$ with bounds $A,B>0$ and pre-frame operator $T$. Then $F$ is $\mu$-complete and the mapping
\begin{align*}
V:\; L^{2} & (\Omega ,A)\longrightarrow L^{2}(\Omega ,A)\\
& \varphi\longmapsto\int_{\Omega}\varphi(\omega)\langle F(\omega),F(.)\rangle d\mu(\omega)
\end{align*}
defines a bounded, adjointable and invertible operator.
\end{corollary}
\begin{proof}
Let $F$ be a continuous Riesz basis for $U$ with bounds $A,B>0$. Then the synthesis operator $T$ satisfies $\Vert T\Vert\leq\sqrt{B}$.
Also,
\begin{align*}
(T^{*}T)(\varphi) & =T^{*}(\int_{\Omega}\varphi(\omega) F(\omega) d\mu(\omega))\\
& = \lbrace\langle\int_{\Omega}\varphi(\omega) F(\omega) d\mu(\omega),F(\gamma)\rangle\rbrace_{\gamma\in\Omega}\\
& = \lbrace\int_{\Omega}\varphi(\omega)\langle F(\omega) ,F(\gamma)\rangle d\mu(\omega)\rbrace_{\gamma\in\Omega}.
\end{align*}
Then $V=T^{*}T$. Moreover, $F$ is Riesz-type and $T^{*}$ is onto. Then by lemma \ref{ISIN}, $V$ is adjointable and invertible operator and 
\begin{equation*}
\Vert (T^{*}T)^{-1}\Vert^{-1}\leq V\leq\Vert T^{*}\Vert^{2}\leq B.
\end{equation*}
\end{proof}    
 According to the Theorem \ref{R-type}, in the next corollary we show the relation between two continuous Riesz bases for a Hilbert $C^{\ast}$-module $U$. 
\begin{corollary}
Let $F,G:\Omega \to U$  be two continuous Riesz bases for Hilbert $C^{\ast}$-module $U$ over a unital $C^*$-algebra $\mathcal A$ and $T_{F},T_{G},S_{G}$ are the pre-frame operator of $F$, the pre-frame operator of $G$ and the frame operator of $G$ respectively. Then there exists an invertible operator $K\in End^{*}(U)$ such that $G=S_{G}K^{*}F$.
\end{corollary}     
\begin{proof}
Let $f\in U$ such that $(T_{G}T^{*}_{F})f=0$. Then $T_{G}((T^{*}_{F}f)(\omega))=0$ for all $\omega\in\Omega$ and $\int_{\Omega}\langle f,F(\omega)\rangle G(\omega) d\mu(\omega)=0$.

Since $G$ is $L^{2}$-independence, so $\langle f,F(\omega)\rangle =0$ for all $\omega\in\Omega$ and by completeness of $F$, we have $f=0$. This shows that $T_{G}T^{*}_{F}$ is injective. Also $F$ a is Riesz-type frame, so $T^{*}_{F}$ is onto. Since $G$ is a continuous frame so $T_{G}$ is onto. Hence $T_{G}T^{*}_{F}$ is onto and so is invertible.

Put $K:=(T_{G}T^{*}_{F})^{-1}$. Then for any $f,g\in U$,
\begin{align*}
\langle f,g\rangle & =\langle K^{-1}Kf,g\rangle\\
& = \langle T^{*}_{F}Kf,T^{*}_{G}g\rangle\\
& = \int_{\Omega}\langle Kf,F(\omega)\rangle \langle G(\omega),g\rangle d\mu(\omega)\\
& = \int_{\Omega}\langle f,K^{*}F(\omega)\rangle \langle G(\omega),g\rangle d\mu(\omega).
\end{align*}
Thus $K^{*}F$ is a dual of $G$. But $G$ is a Riesz-type frame, then $S^{-1}_{G}G=K^{*}F$ and hence $G=S_{G}K^{*}F$.
\end{proof}      

\end{document}